\documentclass[12pt]{article}
\usepackage{amsthm,amsmath,amssymb}
\usepackage{graphicx}
\usepackage[colorlinks=true,citecolor=black,linkcolor=black,urlcolor=blue]{hyperref}
\usepackage[latin1]{inputenc}
\usepackage{subfigure}
\usepackage{epsf, epsfig, tikz, hyperref}
\usepackage{bm}
\usepackage{xifthen}
\usepackage{color}
\tikzstyle{peers}=[circle, draw, fill=black!20, inner sep=0pt, minimum width=8pt]

\newcommand{\ee}[1]{\begin{align} #1 \end{align}}

\newcommand\catalannumber[3]{
  (#1)
  \foreach \dir in {#3}{
    \ifnum\dir=0
    -- ++(1,0)
    \else
    -- ++(0,1)
    \fi
  } |- (#1);
  \draw[dashed] (#1) -- +(2*#2,0);
  \coordinate (prev) at (#1);
  \foreach \dir in {#3}{
    \ifnum\dir=0
    \coordinate (dep) at (1,1);
    \else
    \coordinate (dep) at (1,-1);
    \fi
    \draw[line width=1.5pt,-stealth] (prev) -- ++(dep) coordinate (prev);
  };
}

\theoremstyle{plain}

\theoremstyle{definition}

\theoremstyle{remark}

\title{\bf Counting Restricted Dyck Paths Through Random Walks}

\author{EJ Infeld\\
\small Department of Mathematics\\[-0.8ex]
\small Dartmouth College\\[-0.8ex] 
\small New Hampshire, U.S.A.\\
\small\tt ewa.j.infeld.gr@dartmouth.edu\\
}

\date{}

\begin{document}

\maketitle

\begin{abstract}
We show connection between Dyck paths with peaks of bounded height and random walks. The correspondence between a certain class of random walks and such Dyck paths allows us to develop a probabilistic perspective on Chebyshev polynomials.

  \bigskip\noindent \textbf{Keywords:} Dyck paths, random walks
\end{abstract}

\section{Introduction}
A \textit{Dyck path} of order $k$ is a lattice path with end points $(0, 0)$ and $(2k,0)$, that consists of two kinds of steps: a step up, $(1,1)$ and a step down $(1,-1)$. A \textit{peak} of a Dyck path is the node formed at the joint of an up step and a down step that follows it.

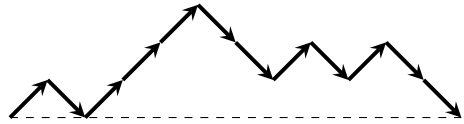
\begin{figure}[ht!]
\centering
\begin{tikzpicture}[scale=0.5]
  \catalannumber{0,0}{6}{0,1,0,0,0,1,1,0,1,0,1,1};
\end{tikzpicture}
\caption{A Dyck path of order 6 with peaks of height 1, 3, 2, and 2 respectively.}
\end{figure}

The number of possible Dyck paths of order $k$ is described by the $k$th Catalan number, $C_k=\frac{1}{k+1}{{2k} \choose k}$. The number $A(n,k)$ of possible Dyck paths of order $k$ that have peak height bounded by $n$ is known to be described by the following generating function: \cite{sets}

\ee{\underset{k\geq 0}{\sum} A(n,k)x^k=\frac{1}{\Large 1+z-\frac{z}{1+z-\frac{z}{\frac{\ddots}{1+z-\frac{z}{1+z-C(z)}}}}}\normalsize,} where the continued fraction has $n+1$ layers and $C(z)$ is the Catalan generating function. Note that the bound on the peak height is equivalent to the condition that in no initial segment the number of steps to the right minus the number of steps up exceeds $n$.

We will present a new and surprising way to count bounded Dyck paths through a property of random walks. In section 2 we discuss the random walk setup and probabilistic results that serve as the starting point. Section 3 is the overview of a calculation that takes us from the realm of probabilities to a suitable generating function. Section 4 presents how the resulting formula outputs the correct functions for some example values of $n$. We present our conclusions in section 5.

\section{The Method}
Consider a path of $m+1$ nodes marked 0 to $m$, and a random walk that starts at node $m-1$ and at every time increment moves one node to the right with probability $0\leq p\leq 1$, $p\not=1/2$ or one node to the left with probability $1-p$: 

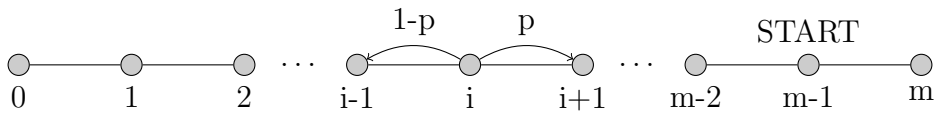
\begin{figure}[ht!]
\centering
\begin{tikzpicture}[scale=0.75]
\foreach \place/\name/\label in {{(-6,0)/a/0}, {(-4, 0)/b/1}, {(-2, 0)/c/2}, {(6,0)/d/m-2}, {(8,0)/e/m-1}, {(10,0)/f/m}, {(2,0)/k/i}, {(4,0)/l/i+1}, {(0,0)/m/i-1}}
    \node[peers,label=below:\label] (\name) at \place{} ;
    \foreach \source/\sink in {a/b, b/c, d/e, e/f, k/l, k/m}
      \path (\source) edge (\sink);
\path[->] (k) edge [bend left]  node[above]{p} (l);
\path[->] (k) edge [bend right]  node[above]{1-p} (m);
\draw[text=black] (-1,0) node {\dots};
\draw[text=black] (5,0) node {\dots};
\draw (8,0.6) node {START};
\end{tikzpicture}
\end{figure}

Suppose that the random walk continues until it reaches either 0 or $m$. Let $L_m(p)$ be the expected hitting time from $m-1$ to $m$, provided the walk does not reach 0. This value can be obtained in a straightforward manner by averaging over the length of all walks that start at $m-1$ and end at $m$ without reaching 0, weighted by the conditional probability of such walk occurring provided that a walk does not in fact hit 0. Since any walk from $m-1$ to $m$ is of odd length, making exactly one more step to the right than to the left, the probability of any such particular walk of length $2k+1$ occurring is $p^{k+1}(1-p)^k$. Let $\Pi_{m,p}$ be the probability that a walk starting at $m-1$ will reach $m$ before reaching 0. Then the conditional probability given that the walk does not reach 0, of any particular walk from $m-1$ to $m$ of length $2k+1$, not reaching 0, occurring is $p^{k+1}(1-p)^k/\Pi_{m,p}$.

For an arbitrary random walk that starts at $m-1$, with probability $p$ the first step is to the right. With probability $\Pi_{m,p}-p$ the first step is to the left, but the walk still reaches $m$ without reaching 0. In that case, it takes an expected $L_{m-1}(p)$ steps to come back to $m-1$, and $L_m(p)$ again to finally reach $m$. When conditioning on reaching $m$ before 0 the whole expression is normalized by $\Pi_{m,p}$. So $L_m(p)$ obeys the following relation: \ee{\Pi_{m,p}L_(p)=p+(\Pi_{m,p}-p)(L_{m-1}(p)+L_m(p)).\label{expval}}

If we were to calculate $L_m(p)$, we would like to sum over the possible lengths $2k+1$, taking into account the number $B(m,k)$ of suitable walks of this length and the conditional probability of each walk occurring. We find that: \ee{L_m(p)=\underset{k\geq 0}{\sum}B(m,k)(2k+1)p^{k+1}(1-p)^k/\Pi_{m,p}\\ =p/\Pi_{m,p}\underset{k\geq 0}{\sum}B(m,k)(2k+1)(p(1-p))^k.} The probability $\Pi_{m,p}$ is \cite{DS}: \ee{\Pi_{m,p}=\frac{1-((1-p)/p)^{m-1}}{1-((1-p)/p)^m}\\ =p\frac{(1-p)^{m-1}-p^{m-1}}{(1-p)^m-p^m}.} Note then that the part of each walk that includes nodes $1,\dots,m-1$ is a collection of an equal number of left and right steps, and at no point the number of left steps already taken minus the number of right steps already taken is more than $m-2$. So each such walk of total length $2k+1$ corresponds to a Dyck path of order $k$ and peak height no more than $m-2$.

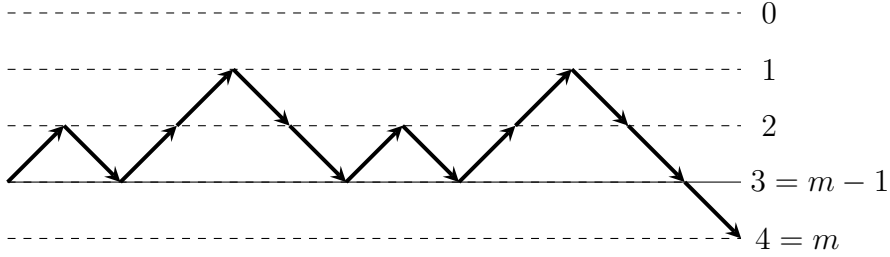
\begin{figure}[ht!]
\centering
\begin{tikzpicture}[scale=0.75]
  \catalannumber{0,0}{6}{0,1,0,0,1,1,0,1,0,0,1,1,1};
  \draw (0,0) -- (13,0);
  \draw[dashed] (0,-1) -- (13,-1);
  \draw[dashed] (0,1) -- (13,1);
  \draw[dashed] (0,2) -- (13,2);
  \draw[dashed] (0,3) -- (13,3);
  \draw[text=black] (13.5,2) node {1};
  \draw[text=black] (13.5,3) node {0};
  \draw[text=black] (13.5,1) node {2};
  \draw[text=black] (14.4,0) node {$3=m-1$};
  \draw[text=black] (14,-1) node {$4=m$};
\end{tikzpicture}
\caption{A Dyck path of order 6 bounded in height by 2 corresponds to a walk of length 7 with $m=4$.}
\end{figure}

 If $A(n,k)$ is the number of distinct Dyck paths of order $k$ bounded above by $n$, we find that $A(n,k)=B(n+2,k)$. If we can then find a suitable expression for $L_m(p)$, we could use the following relation to find  the generating function $\underset{k\geq 0}{\sum}(2k+1)A(n,k)x^k$: \ee{\underset{k\geq 0}{\sum}A(n,k)(2k+1)(p(1-p))^k=\Pi_{n+2,p}/p\times L_{n+2}(p),} where:\ee{\Pi_{m,p}=p\frac{(1-p)^{m-1}-p^{m-1}}{(1-p)^m-p^m}.}
 
\section{Calculation}
We will now focus on obtaining a suitable form for $\Pi_{n+2,p}/p\times L_{n+2}(p)$ as a function of $x=p(1-p)$.  We will use the following notation: 
\begin{itemize}
\item $g_i(p)=(1-p)^i-p^i$,
\item $f_i(p)=\Pi_{i,p}/p=\frac{(1-p)^{i-1}-p^{i-1}}{(1-p)^i-p^i}=\frac{g_{i-1}(p)}{g_i(p)}$,
\item $h_i(p)=\frac{g_i(p)}{g_1(p)}=\frac{g_i(p)}{1-2p}$,
\end{itemize} It is worth noting that: \ee{g_i(p)=g_i(p)(1-p+p)=g_{i+1}(p)+xg_{i-1}(p),} and that equation (\ref{expval}) is equivalent to: \ee{L_m(p)=\Pi_{m,p}/p+L_{m-1}(p)(\Pi_{m,p}/p-1)\\ =f_m(p)+L_{m-1}(p)(f_m(p)-1).} By iterating the above relation we can discover that: $$L_m(p)=f_m(p)+f_{m-1}(p)(f_m(p)-1)+$$ $$+f_{m-2}(p)(f_{m-1}(p)-1)(f_{m-2}(p)-1)+\dots+L_2(p)\underset{j=3}{\overset{m}{\prod}}(f_j(p)-1),$$ where $L_2(p)=1$. We can express this in a more transparent form as:
$$L_m(p)=\underset{i=2}{\overset{m}{\sum}}f_i(p)\underset{j=i+1}{\overset{m}{\prod}}(f_j(p)-1),$$ since $f_2(p)=L_2(p)=1.$ Finally, since we have established that: \ee{\underset{k=o}{\overset{\infty}{\sum}}(p(1-p))^k(2k+1)A(m-2,k)=f_m(p)L_m(p),} we would like to express $f_m(p)L_m(p)$ as a formal power series of $p(1-p)=x$.
\ee{f_m(p)L_m(p)=\frac{g_{m-1}(p)}{g_{m}(p)}\left[\underset{i=2}{\overset{m}{\sum}}\frac{g_{i-1}(p)}{g_{i}(p)}\underset{j=i+1}{\overset{m}{\prod}}(\frac{g_{j-1}(p)}{g_{j}(p)}-1)\right]\\ =\frac{g_{m-1}(p)}{g_{m}(p)}\left[\underset{i=2}{\overset{m}{\sum}}\frac{g_{i-1}(p)}{g_{i}(p)}\underset{j=i+1}{\overset{m}{\prod}}x\frac{g_{j-2}(p)}{g_{j}(p)}\right],} since $xg_{j-2}(p)=g_{j-1}(p)-g_{j}(p)$. This can be simplified to:\ee{f_m(p)L_m(p)=\frac{1}{(g_{m}(p))^2}\underset{i=2}{\overset{m}{\sum}}x^{m-i}(g_{i-1}(p))^2=\frac{1}{(g_{m}(p))^2}\underset{i=1}{\overset{m-1}{\sum}}x^{m-i-1}(g_{i}(p))^2.}

The rest of the calculation is as follows: \ee{g_i(p)^2=(1-p)^{2i}+p^{2i}-x^i=x^i((\frac{1-p}{p})^i)+(\frac{p}{1-p})^i-2),} \ee{\underset{i=1}{\overset{m-1}{\sum}}x^{m-i-1}(g_{i}(p))^2=\underset{i=1}{\overset{m-1}{\sum}}x^{m-1}\left[(\frac{1-p}{p})^i+(\frac{p}{1-p})^i-2\right]\\ =x^{m-1}\left[\frac{1-p}{p}\frac{1-(\frac{1-p}{p})^m-1}{1-\frac{1-p}{p}}+\frac{p}{1-p}\frac{1-(\frac{p}{1-p})^m-1}{1-\frac{p}{1-p}}-2(m-1)\right]\\ =x^{m-1}\left[\frac{1-p-(1-p)(\frac{1-p}{p})^m-1}{2p-1}+\frac{p-p(\frac{p}{1-p})^m-1}{1-2p}-2(m-1)\right] \\ = \frac{x^{m-1}(1-2m)g_1(p)+(1-p)^{2m-1}-p^{2m-1}}{g_1(p)}=\frac{x^{m-1}(1-2m)g_1(p)+g_{2m-1}(p)}{g_1(p)}.} We can therefore conclude that: \ee{f_m(p)L_m(p)=\frac{x^{m-1}(1-2m)g_1(p)+g_{2m-1}(p)}{g_{m}(p)^2g_1(p)}\\ =\frac{x^{m-1}(1-2m)/g_1(p)^2+(g_{2m-1}(p)/g_1(p))/g_1(p)^2}{(g_{m}(p)/g_1(p))^2}\\ =\frac{x^{m-1}(1-2m)+h_{2m-1}(p)}{g_1(p)^2h_m(p)^2}.} Finally, note that $g_1(p)=(1-2p)^2=1-4p+4p^2=1-4x$, and that since $h_1(p)=h_2(p)=1$, and $h_m(p)=h_{m-1}(p)-xh_{m-2}(p)$ for $m\geq 3$, we can define $H(x)=h(p)$ as: \begin{itemize}\item $H_1(x)=H_2(x)=1$ \item For $m\geq 3$, $H_m(x)=H_{m-1}(x)-xH_{m-2}(x),$\end{itemize} and: \ee{f_m(p)L_m(p)=\frac{x^{m-1}(1-2m)+H_{2m-1}(x)}{(1-4x)H_m(x)^2}.} Therefore, for $A(n,k)$ the number of Dyck paths of order $k$ and height no more than $n$:
\ee{\underset{k\geq 0}{\sum}A(n,k)(2k+1)x^k=\frac{x^{n+1}(-2n-3)+H_{2n+1}(x)}{(1-4x)H_{n+2}(x)^2}.\label{gf}}

This result is equivalent to the continued fraction form that is known, but the probability viewpoint is novel. The factor of $2k+1$ could be removed by integrating both sides, but the formula on the right would lose much of its simplicity. We prioritize efficient use of computational resources and so propose that the factor should be left as is.

\section{Examples}
Let us check this expression for some small numbers. Note that $f_2(p)L_2(p)=0+1=1$, and in fact the only Dyck path of height 0 is the empty path. 

\subsection{Paths of height 1}
For any positive length, there is only one path of height 1, and so we are hoping to obtain: \ee{f_3(p)L_3(p)=\underset{k\geq 0}{\sum}(2k+1)x^k.\label{path3}} We do it as follows. From \ee{f_3(p)L_3(p)=\frac{x^{2}(-5)+H_{5}(x)}{(1-4x)H_3(x)^2}}  we find that $H_3(x)^2=(1-x)^2$, and $H_5(x)=1-3x+x^2$. 

Since $H_5(x)-5x^2=1-3x-4x^2=(1-4x)(1+x)$, we conclude that: \ee{f_3(p)L_3(p)=\frac{1+x}{(1-x)^2}=\underset{k\geq 0}{\sum}(2k+1)x^k.}

\subsection{Paths of height at most 2}
Similarly, we can check this formula for Dyck paths of heights no more than 2 as follows: \ee{f_4(p)L_4(p)=\frac{x^{2}(-7)+H_{7}(x)}{(1-4x)H_4(x)^2}=\frac{1-5x+6x^2-8x^3 }{(1-4x)(1-2x)^2 }=\frac{1-x+2x^2}{(1-2x)^2 }\\ =(1-x+2x^2)\underset{k\geq 0}{\sum}(k+1)(2x)^k=1+\underset{k\geq 1}{\sum}(2k+1)2^{k-1}x^k.}

To see that this is the generating function we are looking for, consider such path of some order $k\geq 1$. This path can touch 'ground' after any even number of steps, and necessarily touches ground after 0 and after $2k$ steps. There are $k-1$ spots where it could touch ground or not. And in fact for any subset of these, there is a unique Dyck path of length $2k$ of height at most 2. Therefore there are $2^{k-1}$ such paths of order $k$.

\section{Conclusions}

We found a generating function that describes the number of Dyck paths of order $k$ restricted in height by $n$, which is (\ref{gf}). An iterative method to count these already existed, but we believe that the random walk connection is of interest.

\end{document}